\documentclass{article}

\begin{document}

Tsemo Aristide

5155 Avenue de Gaspe, Appt 214

H2T 2A1

Quebec, Canada

 tsemoaristide@hotmail.com

\bigskip

\bigskip

{\bf Tours de torseurs g\'eometrie differentielle des suites de
compositions de fibr\'es principaux, et th\'eorie des cordes.}

\bigskip

\bigskip

\centerline{\bf R\'esum\'e.}

Soient $f:M\rightarrow N$ un fibr\'e principal de groupe
structural $L$, et $1\rightarrow H_1\rightarrow L_1\rightarrow
L\rightarrow 1$ une extension de $L$ \`a noyau commutatif $H_1$.
Le probl\`eme d'\'etendre le groupe structural de $f$, en $L_1$,
est formul\'e par la th\'eorie des gerbes. Nous \'etudions la
situation plus g\'en\'erale suivante:  consid\'erons une famille
d'extensions \`a noyaux commutatifs $1\rightarrow
H_{i+1}\rightarrow L_{i+1}\rightarrow L_i\rightarrow 1$, $0<i<n$.
Le probl\`eme de $n-$rel\`evement est celui d'\'etendre le groupe
structural de $f$ \`a $L_n$. Il chouk\`ele la notion de tour de
torseurs que nous d\'efinissons, et \'etudions la g\'eom\'etrie
diff\'erentielle. Nous d\'efinissons notamment les notions de
structures connective, et holonomie, qui permettent de
g\'en\'eraliser l'action de Wess Zumino Witten (WZW) aux branes.
La th\'eorie naturelle plus g\'en\'erale qui permettrait de
formuler ce probl\`eme est celle des $n-$cat\'egories, les
nombreuses applications d'une telle th\'eorie sont mentionn\'ees
dans les travaux de Alexandre Grothendieck, malheureusement elle
est incomprise actuellement.

\bigskip

\bigskip

\centerline{\bf Introduction.}

\bigskip
\bigskip

Soient $N$ une vari\'et\'e differentielle, $L$ un groupe de Lie,
et $f:M\rightarrow N$ un fibr\'e principal de fibre type $L$.
Consid\'erons la suite exacte:

$$
1\longrightarrow H\rightarrow L_1\rightarrow L\rightarrow 1
$$

o\`u $H$ est un groupe commutatif. Le probleme initial \'etudi\'e
appel\'e encore probl\`eme de $1-$rel\`evement, est celui de
l'existence  d'un fibr\'e  principal $M_1$, de fibre type $L_1$,
et de base $N$, tel que le quotient de $M_1$ par $H$ est $M$.
L'existence d'un tel fibr\'e n'est pas toujours assur\'ee.
L'obstruction \`a son existence est donn\'ee par le cocycle
classifiant d'une gerbe dont l'\'etude de la g\'eom\'etrie
diff\'erentielle est bien connue.

Consid\'erons maintenant une succession de suites exactes de
groupes de Lie

$$
0\rightarrow H_1\rightarrow L_1\rightarrow L\rightarrow 0
$$
...
$$
0\rightarrow H_{i+1}\rightarrow L_{i+1}\rightarrow L_i\rightarrow
0
$$
...
$$
0\rightarrow H_n\longrightarrow L_{n}\longrightarrow
L_{n-1}\rightarrow 0
$$
telles que les $H_i$ soient  des groupes commutatifs. Le
probl\`eme de $n-$rel\`evement est de d\'eterminer les
obstructions, et la g\'eom\'etrie diff\'erentielle du probl\`eme
de l'existence de
  suites de fibr\'es principaux

$$
f_{i}:M_{i}\rightarrow M_0=N, i=0,..n.
$$

telle que l'application $f_i$ est une fibration principale de
fibre type $L_{i}$, et le quotient de $M_{i}$ par $H_i$ est
$M_{i-1}$.

 L'existence des applications $f_i$ entraine celle du diagramme commutatif

$$
\matrix{M_{n}{\buildrel {f_n}\over{\longrightarrow}}
M_{n-1}...{\buildrel {f_n}\over{\longrightarrow}}N\cr\downarrow
&&\downarrow \cr N{\buildrel {Id}\over{\longrightarrow}}N...
{\buildrel {Id}\over{\longrightarrow}}N}
$$

Les th\'eories de faisceaux de cat\'egories, et de gerbes, ont
\'et\'e d\'efinies par J. Giraud, et utilis\'ees pour r\'esoudre
le premier probl\`eme de rel\`evement ont pour objectif de fournir
le cadre th\'eorique  aux probl\`emes de recollement de structures
d'ordre $1$.

Deux m\'ethodes fondamentales peuvent \^etre mises en exerguent
dans la construction de structures en math\'ematiques ce sont:

Le processus de compl\'etion, qui consiste \`a l'extension de
structures, et d'espaces o\`u des probl\`emes non r\'esolus
trouveront solutions. Ce type de construction caract\'erise le
passage de l'alg\`ebre \`a l'analyse.

\medskip

La seconde proc\'edure est celle du recollement utilis\'ee dans le
pr\'esent probl\`eme. Pour r\'esoudre les probl\`emes de
recollement, ou de mani\`ere \'equivalente, pour recoller des
objets d\'efinis localement, doivent \^etre d\'efinies des
th\'eories cohomologiques interpr\'et\'ees g\'eom\'etriquement.

Cette derni\`ere  m\'ethode de construction s'illustre dans
l'exemple suivant:

Une vari\'et\'e topologique $N$ est construite, en recollant des
ouverts d'un espace vectoriel $E$ \`a l'aide d'homeomorphismes.
Pour \'etudier $N$, on g\'en\'eralise les outils d\'efinis dans
l'\'etude de l'espace vectoriel $E$. C'est ainsi que les fonctions
continues d\'efinies localement peuvent \^etre d\'efinies sur $N$.
Ce probl\`eme de $1-$recollement, est r\'esolu par
l'interpr\'etation g\'eom\'etrique des groupes de cohomologies
d'ordre $1$ de la notion de faisceau invent\'e par J. Leray. La
notion de faisceau conduit \`a celle de fibr\'es principaux
localement triviaux. Est-il possible de recoller des fibr\'es
d\'efinis localement? Pour r\'esoudre ce probl\`eme de
$2-$recollement, on utilise les notions de faisceaux de
cat\'egories et de gerbes d\'efinies par J. Giraud.

Plus g\'en\'eralement le probl\`eme de recollement  des faisceaux
de cat\'egories conduit \`a celui de d\'efinir la notion de
faisceaux de $3-$cat\'egories,..., celui  de recollement des
faisceaux de $n-$cat\'egories conduit \`a celui de d\'efinir la
notion de faisceaux de $n+1-$cat\'egories. Malheuresement il
n'existe pas de th\'eories satisfaisantes de $n-$cat\'egories
$n>2$, bien qu'une telle th\'eorie ait \'et\'e initi\'ee dans
Tsemo[T2]. On contourne cette difficult\'e en introduisant la
notion de tour de torseurs qui permet de donner une
interpr\'etation g\'eom\'etrique des classes de cohomologies
ab\'eliennes. Une tour de torseurs est d\'efinie par une suite de
foncteurs $F_n\rightarrow F_{n-1}...\rightarrow F_1$, \`a laquelle
est associ\'ee une famille de faisceaux $L_2$,...,$L_n$ d\'efinis
sur $F_1$ qu'on suppose munie d'une topologie. Ces donn\'ees
doivent satisfaire des contraintes refl\'etant la notion de
recollement. A une tour de torseurs, est associ\'ee un ensemble de
classes de cohomologies $[c_2]\in H^2(F_1,L_2)$,..., $[c_{n}]\in
H^{n}(F_1,L_n)$ d\'eterminant la classe d'isomorphisme de la tour.

On associe au probl\`eme de $n-$rel\`evement une tour de torseurs
$F_n\rightarrow F_{n-1}...F_2\rightarrow F_1$,   dont on \'etudie
la g\'eom\'etrie diff\'erentielle.  Les objets des cat\'egories
$F_i$ sont des fibr\'es principaux au-dessus d'ouverts de $N$. On
\'etend  la notion de structure connective introduite par
Brylinski \`a ce cadre, et on d\'efinit  la courbure d'une
structure connective, ainsi que les classes charact\'eristiques. A
cette tour de torseurs, est associ\'ee des formes
diff\'erentielles de $N$, d\'efinies \`a l'aide  d'un recouvrement
$(U_i)_{i\in I}$ de $N$ par des ouverts simplement contractibles:

\medskip

 La courbure $L$ est une $n+1-$forme.  Pour tout ouvert
$U_{i_1...i_l}$, on d\'efinit une $n-l+1$ forme $u_{i_1..i_l}$. La
relation suivante est v\'erifi\'ee:

$$
d(u_{i_1..i_l})=\sum_{j=1}^{j=l}(-1)^ju_{i_1..\hat i_j..i_l}
$$
 pour $l<n+1$, et

 $$
 \sum_{j=1}^{j=n+1}(-1)^ju_{i_1..\hat i_j..i_{n+1}}=(c_{i_1..i_{n+1}})^{-1}dc_{i_1..i_{n+1}}
$$

En physique th\'eorique moderne, la particule a \'et\'e
remplac\'ee par la notion de corde. L'action d\'eterminant cette
\'evolution dans le mod\`ele de WZW est d\'efinie par
 l'holonomie d'une structure connective d'une gerbe. Lorsque la corde
 est soumise \`a des contraintes ou branes, la formulation variationnelle du mouvement
 de la brane  est d\'efinie par l'\'evolution d'une
sous-vari\'et\'e $N'$ de $N$ sous une action  d\'etermin\'ee par
un ensemble de formes diff\'erentielles satisfaisant les
propri\'et\'es mentionn\'ees ci-dessus. Nous d\'efinissons
l'holonomie d'une structure connective d'une tour de torseurs.
Ceci permet de g\'en\'eraliser l'action de WZW en dimenension
sup\'erieure.

\medskip

Les probl\`emes de recollements d'objets tirent leur origine dans
les travaux d'Alexander Grothendieck. Andr\'e Weil avait
\'enonc\'e les conjectures de Weil concernant la fonction zeta, et
avait sugg\'er\'e que leurs d\'emonstrations r\'esulteraient de
formule de Leftchetz d'une th\'eorie de cohomologie ad\'equate.
Jean Pierre Serre avait remarqu\'e que les fibr\'e principaux
au-dessus des vari\'et\'es alg\'ebriques \'etait localement
triviaux \`a un rev\^etement \'etale pr\`es. Motiv\'e par ces
id\'ees, Alexander Grothendieck a d\'efini une notion de
cohomologie pour les cat\'egories qu'il a appliqu\'e ensuite \`a
la g\'eom\'etrie alg\'ebrique cela a donn\'e naissance \`a la
cohomologie \'etale qui a permit a Pierre Deligne de d\'emontr\'e
les conjectures de Weil. Les succ\`es de cette approche ont
sugg\'er\'e \`a Alexandre Grothendieck de g\'en\'eraliser les
th\'eories d'homotopies au vari\'et\'es alg\'ebriques. Une telle
g\'en\'eralisation n\'ec\'essite une d\'efinition de notion de
$n-$groupoides et $n-$gerbes dont de $n-$cat\'egories, la notion
de lacet n'existant pas en g\'eom\'etrie alg\'ebrique. En rang $1$
tout de moins, le groupe fondamemtal se d\'efinit comme objet
universel dans la cat\'egorie des faisceaux localement constant.
Grothendieck ceci a permi a Grothendieck de d\'efinir le groupe
fondamental en topologie \'etale, obtenant une suite exacte ayant
de nombreuses applications en arithm\'etiques. Voedvodsky a
d\'efini une homotopie motivique qui lui a permi de d\'emontr\'e
une conjecture de Milnor.

\bigskip

\bigskip

{\bf I. Etude alg\'ebrique.}

\bigskip

L'objectif de cette partie est de d\'efinir la notion de tour de
torseurs ou tour tors\'ee, et d'en d\'eduire les principales
propri\'et\'es alg\'ebriques. On commence d'abord par rappeller
celle de gerbe, introduite par J. Giraud motiv\'e par le
probl\`eme de $1-$relevement. La notion de gerbe, se d\'efinit
dans le cadre g\'en\'eral de topologies des cat\'egories que nous
pr\'esentons maintenant.

\medskip

{\bf Definitions 1.}

Soit $C$ une cat\'egorie:

 Un crible $U$ est une sous-classe de la classe des objects
$Obj(C)$ telle que pour tout \'el\'ement $u$ de $U$, et toute
fl\`eche $u'\rightarrow u$,  $u'$ est un \'el\'ement de $U$.

\medskip

Une topologie sur $C$ est d\'efinie par la donn\'ee pour tout
objet $u$ de $C$, d'une famille de cribles $J(u)$, dont les
\'el\'ements sont appel\'es raffinements de $u$, v\'erifiant les
propri\'et\'es suivantes:

Pour toute fl\`eche $f:u'\rightarrow u$, et tout crible $U$ de
$J(u)$, $U^f=\{u": u"\rightarrow u'\rightarrow u\in J(u)\}$ est un
\'el\'ement de $J(u')$.

Soit $U$ un crible de $u$, si pour toute fl\`eche $u'\rightarrow
u$, $U^f\in J(u')$, alors $U$ est un \'el\'ement de $J(u)$.

\medskip

L'exemple fondamental de topologies dans une cat\'egorie est le
suivant: Soit $N$ un espace topologique. La cat\'egorie de ses
ouverts $Ouv(N)$  est munie de la topologie dont les raffinements
d'un ouvert $U$, sont les ensembles d'ouverts $(U_i)_{i\in I}$,
tels que $\bigcup_{i\in I}U_i=U$.

\bigskip

{\bf D\'efinition 2.}

Soit $(C,J)$ une cat\'egorie $C$, munie d'une topologie $J$. Un
faisceau d\'efini sur $C$, est un foncteur contravariant de $C$
dans la cat\'egorie des ensembles $Set$

$$
F:C\rightarrow Set,
$$

 tel que pour tout objet $u$, et tout crible $U$ de $J(u)$,

$$
F(u)=lim_{u'\in U} F(u').
$$

\medskip

Nous allons \`a pr\'esent rappeler comment les propri\'et\'es de
recollement
 d' objets  s'expriment par la notion de faisceaux de cat\'egories.

\medskip

{\bf Definitions 3.}

- Soient $f:F\rightarrow E$ un foncteur, $x,y$ deux objets de $F$,
et $n:x\rightarrow y$ une fl\`eche dont l'image par $f$ est la
fl\`eche $N:X\rightarrow Y$ de $E$. On dira que la fl\`eche $n$
est cart\'esienne si et seulement si, pour tout object $z$ de la
fibre au-dessus de $X$, c'est \`a dire tel que $f(z)=X$,

$$
Hom_X(z,x)\rightarrow Hom_N(z,y)
$$

$$
m\rightarrow nm
$$
est bijective, o\`u $Hom_Z(z,x)$ est l'ensemble des fl\`eches
$u:z\rightarrow x$ telles que $f(u)=Id_X$, et $Hom_N(z,y)$ est
l'ensemble des fl\`eches $u':z\rightarrow y$ telles que $f(u')=N$.

\medskip

- On dira que la cat\'egorie  $F\rightarrow E$ est cart\'esienne,
si pour toute fl\`eche $N:X\rightarrow Y$ de $E$, il existe une
fl\`eche car\'esienne $n:x\rightarrow y$ dont l'image par $f$ est
$N$, et la composition de deux fl\`eches cart\'esiennes est
cart\'esienne.

\medskip

- Soient $f:F\rightarrow E$ et $f':F'\rightarrow E$ deux
cat\'egories cat\'esiennes. Un foncteur  $u:F'\rightarrow F$ est
cart\'esien si et seulement si $u$ respecte les fibres et l'image
d'une fl\`eche cart\'esienne est cart\'esienne. On d\'enote par
$Cart_E(F',F)$ la cat\'egorie des foncteurs cart\'esiens entre
$F'$ et $F$. Un exemple de foncteur cart\'esien est l'immersion
$E_X\rightarrow E$, o\`u $E_X$ est la sous-cat\'egorie des objets
de $E$ au-dessus de $X$.

\medskip

{\bf D\'efinition 4.}

On dira que le foncteur $f:F\rightarrow E$ est un faisceau de
cat\'egories, si et seulement si pour tout objet $X$ de $E$, et
tout crible $U$ de $J(X)$,  le foncteur restriction:

$$
Cart_E(E_X,F)\longrightarrow Cart_E(U,F)
$$
est une \'equivalence de cat\'egories. $(E,J)$ est appel\'ee la
base du faisceau de cat\'egories.

\medskip

{\bf Proposition 1.}

Soit $(E,J)$ un topos muni d'une topologie $J$ engendr\'ee par la
famille $(X_i\rightarrow X)_{i\in I}$, un faisceau de cat\'egories
de base $(E,J)$ est \'equivalent \`a la donn\'ee d'une
correspondance

$$
C:U\rightarrow C(U),
$$
o\`u $U$ d\'esigne un objet de $E$, et $C(U)$ une cat\'egorie. Les
propri\'et\'es suivantes sont satisfaites: Pour toute fl\`eche
$U\rightarrow V$ de $E$, il existe un foncteur

$$
r_{U,V}:C(V)\rightarrow C(U)
$$

telle que

$$
r_{U,V}\circ r_{V,W}\simeq r_{U,V}.
$$

En pratique on supposera toujours que la pr\'ec\'edente
\'equivalence est une \'egalit\'e.

(i) Recollement des objets:

 Soient $X$ un objet de $E$, $(X_i\rightarrow X)_{i\in I}$ un recouvrement
 de $X$,  $x_i$ un objet de $C(X_i)$, et

$$
g_{ij}:r_{X_i\times_XX_j,X_j}(x_j)\rightarrow
r_{X_i\times_XX_j,X_i}(x_i)
$$

une famille de fl\`eches telles que $g_{ij}g_{jk}=g_{ik}$, il
existe un objet $x$ de $C(X)$ tel que $x_i=r_{X_i,X}(x)$.

(ii) Recollement des fl\`eches:

Pour tous objets $x$ et $y$ de $C(U)$, le prefaisceau

$$
U\rightarrow Hom(r_{U,X}(x),r_{U,X}(y))
$$
est un faisceau de la cat\'egorie au-dessus de $X$ munie de la
topologie induite.

\medskip

Si  de plus les propri\'et\'es suivantes sont satisfaites

(i') la topologie $J$ est engendr\'ee par une famille $(X_i)_{i\in
I}$, telle que $C(X_i)$ n'est pas vide.

(ii') Pour tout objets $x$ et $y$,  de $C(X)$, $r_{X_i,X}(x)$ et
$r_{X_i,X}(y)$ sont isomorphes.

(iii') Il existe un faisceau $L$ d\'efini sur $C$, tel que pour
tout objet $x$ de $C(X)$, $Hom(x,x)=L(X)$,

alors le faisceau de cat\'egories est appel\'e gerbe, et $L$ son
lien.

\bigskip

{\bf Le cocycle associ\'e \`a une gerbe.}

\medskip

D\'efinissons maintenant le cocycle classifiant associ\'e \`a une
gerbe dont le lien est un faisceau en groupe ab\'elien.

 Pour tout \'el\'ement $X_i$ de la famille g\'en\'eratrice $(X_i)_{i\in
I}$, on consid\`ere un objet $x_i$ de $C(X_i)$, et une fl\`eche

$$
g_{ij}:r_{X_i\times_XX_j,X_j}(x_j)\rightarrow
r_{X_i\times_XX_j,X_i}(x_i).
$$

 On d\'efinit alors

$$
c_{ijk}=g_{ki}g_{ij}g_{jk}.
$$

\bigskip

{\bf Th\'eoreme 1. Giraud [Gi] p. 264.}

{\it La famille $c_{ijk}$ est un $2-$cocycle de Cech. L'ensemble
des gerbes sur $E$ de lien $L$ est isomorphe \`a $H^2(E,L)$.}

\bigskip

 Rappelons maintenant comment cette th\'eorie s'applique
   au probl\`eme de $1-$rel\`evement des fibr\'es principaux.

 Consid\'erons un fibr\'e principal localement trivial
$f:M\rightarrow N$ de groupe structural $L$, et une extension
$1\rightarrow H\rightarrow L_1\rightarrow L\rightarrow 1$ \`a
noyau ab\'elien $H$. On supppose que le morphisme $L_1\rightarrow
L$ a des sections locales.

 D\'efinissons la gerbe $C$ sur la cat\'egorie des ouverts de $N$,
qui a tout ouvert $U$ de $N$ associe la cat\'egorie $C(U)$, dont
les objets sont les fibr\'es principaux localement triviaux de
groupe structural $L_1$ au-dessus de $U$, tels que le quotient de
tout objet $e$ de $C(U)$ par $H$ est la restriction de $f$ \` a
$U$. Les morphismes entre \'el\'ements de $C(U)$ sont les
morphismes de fibr\'es principaux qui se projettent sur
l'identit\'e de $f_{\mid U}$.

\medskip

Le cocycle associ\'e \`a la gerbe $C$ est \`a valeurs dans le
faisceau des sections du $H-$fibr\'e principal d\'efini sur $N$
par $f$. Il se d\'efinit comme suit:

On consid\`ere un recouvrement ouvert $(U_i)_{i\in I}$ par des
ouverts contractibles. On choisit dans toute cat\'egorie $C(U_i)$
un objet $x_i$, et un morphisme $g_{ij}:x_j^i\rightarrow x_i^j$.
Le cocycle est alors d\'efini par $c_{ijk}=g_{ki}g_{ij}g_{jk}$.

\bigskip

{\bf 2. Le cas g\'en\'eral.}

\medskip

Dans un de nos travaux r\'ecents [T2], nous avons d\'efini la
notion de tour de gerbes pour r\'esoudre les probl\`emes
d'extensions. Cette notion n'est pas utile dans toute sa
g\'en\'eralit\'e pour \'etudier le probl\`eme de $n-$rel\`evement
des fibr\'es principaux. Dans la suite on va d\'evelopper une
notion moins g\'en\'erale qu'on appelle tour de torseurs ou tour
tors\'ee.

\medskip

{\bf D\'efinition 5.}

Une tour de torseurs est une suite de foncteurs

$$
F_n{\longrightarrow}F_{n-1}...F_2{\longrightarrow} F_1,
$$
Soit $U$ un objet de $F_1$, on note ${F_p}_U$, les objets de $F_p$
se projettant sur $U$ par la suite de foncteurs $F_p\rightarrow
F_1$. Pour toute fl\`eche $U\rightarrow V$, il existe un foncteur
restriction ${r^p}_{U,V}:{F_p}_V\rightarrow {F_p}_U$ tel que
${r^p}_{U_1,U_2}\circ {r^p}_{U_2,U_3}={r^p}_{U_1,U_3}$ v\'erifiant
les propri\'et\'es suivantes:

(1) Le foncteur $p_2:F_2\rightarrow F_1$ est une gerbe de lien
$L_2$.

(2) Les fl\`eches de la cat\'egorie $F_i$, se projettant sur
l'identit\'e par le foncteur $p_i:F_i\rightarrow F_{i-1}$ $i>1$
sont inversibles.

(3) Il existe des faisceaux $L_2$,...,$L_n$ d\'efinis sur $F_1$
tels que pour tout objet $x_i$ de $F_i$, et $x_{i+1}$ de
${F_{i+1}}_{x_i}$, il existe un isomorphisme:

$$
Aut_{x_i}(x_{i+1})\rightarrow L_{i+1}(p_2..p_{i+1}(x_{i+1}))
$$

 Cet isomorphisme est naturel dans le sens qu'il existe un
isomorphisme:

$$
L_{i+1}(p_2..p_{i+1}(x_{i+1}))\rightarrow
Aut_{x_i}({{F_{i+1}}_{x_i}}).
$$

$$
f\longrightarrow u(f)
$$

 o\`u $Aut_{x_i}({{F_{i+1}}_{x_i}})$ est le groupe des automorphismes de
 ${F_{i+1}}_{x_i}$ qui se projettent sur l'identit\'e. Pour tout
objet  $x_{i+1}$ de ${F_{i+1}}_{x_i}$, l'action de $u(f)$ sur
$x_{i+1}$ est celle d\'efinie par l'isomorphisme pr\'ec\'edent
entre $L_{i+1}(p_2..p_{i+1}(x_{i+1}))$ et $Aut_{x_i}(x_{i+1})$.

(4) Les fl\`eches de $F_i$ se projettant sur l'identit\'e sont
cart\'esiennes. Pour toute fl\`eche $n:x_i\rightarrow y_i$, de la
cat\'egorie $F_i$, se projettant sur l'identit\'e. et la
composition de deux tels relev\'es cart\'esiens est cart\'esienne.
Il existe donc une fl\`eche

$$
n^*:x_{i+1}\rightarrow y_{i+1},
$$

au-dessus de $n$ telle que pour toute fl\`eche $m:y_i\rightarrow
z_i$ de $F_i$,

$$
(mn)^*=u(m,n)m^*n^*
$$

o\`u $u(m,n)$ d\'esigne un automorphisme de $z_{i+1}$ se
projettant sur l'identit\'e.

(5) Connectivit\'e locale

Il existe un recouvrement $(U_i\rightarrow U)_{i\in I}$ de $U$
telle que  Pour tout \'el\'ements $x_{l+1}$ et $y_{l+1}$ de
${F_{l+1}}$, on suppose que
$p_{l+1}...p_2(x_{l+1})=p_{l+1}...p_2(y_{l+1})=U_i$ il existe un
isomorphisme $u_l:x_{l+1}\rightarrow y_{l+1}$ qui se projette sur
l'identit\'e de $p_{l+1}..p_c(x_{l+1})$ si
$p_{l+1}..p_c(x_{l+1})=p_{l+1}..p_c(y_{l+1})$. Cet automorphisme
commute avec l'action de $L_{l+l}$ de $(3)$ et $L'_l$ voir $(6)$.

\medskip

(6) La tour tors\'ee est ab\'elienne si et seulement si les
faisceaux $L_i$ sont commutatifs, $Aut({F_{i+1}}_{x_{i}})$ $i>1$,
le groupe des automorphismes de la fibre
 ${F_{p+1}}_{x_{p-1}}$,  se projettant sur l'identit\'e de
 $x_{p-1}$ est un faisceau en  groupes ab\'eliens $L'_p$ au-dessus de $F_1$ commutant avec les
 isomorphismes d\'efinis en $(5)$.

(7) Il existe une famille couvrante $(X_i\rightarrow X)_{i\in I}$
telle que pour tout $i$, ${F_l}_{X_i}$ n'est pas vide.

 (8) L' application
$n\rightarrow m^*n^*$ est lin\'eaire.

\medskip

{\bf D\'efinition 6.}

Soit $U$ un objet de $F_1$, la fibre de $U$ est la $n-1-$tour de
torseurs $F'_{n-1}\rightarrow...\rightarrow F'_1$,  telle que
$F'_1$ est la fibre ${F_2}_U$, supposons d\'efinie la cat\'egorie
$F'_i$, les objets de $F'_{i+1}$ sont les objets des cat\'egories
${F_{i+2}}_{e'_i}$ o\`u $e'_i$ est un objet de $F'_i$. Les
morphismes et extensions sont
 induits par ceux de la tour $F_n\rightarrow..\rightarrow F_1$.

\medskip

{\bf D\'efinition 7.}

Soient ${T_i=F_n}^i\rightarrow {F_{n-1}}^i\rightarrow..\rightarrow
{F_1}^i$, $i=1,2$ deux tours tors\'ees de m\^eme liens
$(L_2,..,L_n)$ telles que ${F_1}^1$ est ${F_1}^2$. Un
$n-i$morphisme entre $T_1$ et $T_2$, est une famille de foncteurs
$u_n:{F_n}^1\rightarrow
{F_n}^2$,...,$u_{n-i}:{F_{n-i}}^1\rightarrow {F_{n-i}}^2$ tels que
le diagramme suivant soit commutatif:

$$
\matrix{{F_n}^1{\buildrel{f_n}^1\over{\longrightarrow}}
{F_{n-1}}^1...{F_{n-i}}^1{\buildrel{f_{n-i}}^1\over{\longrightarrow}}
{F_{n-i-1}}^1 \cr u_n\downarrow  u_{n-1}\downarrow \downarrow
u_{n-i+1}\cr
{F_n}^2{\buildrel{f_n}^2\over{\longrightarrow}}{F_{n-1}}^2...{F_{n-i}}^2
{\buildrel{f_{n-i}}^1\over{\longrightarrow}}{F_{n-i-1}}^2}
$$

On suppose de plus que le foncteur $u_1$ est continu.

Un $n-1-$morphisme est un isomorphisme si et seulement si il
existe un morphisme $(v_1,...,v_n):T_2\rightarrow T_1$ tel que
$v_iu_i\simeq Id_{F^1_i}$,   et $u_iv_i\simeq Id_{F^2_i}$, les
isomorphismes entre les foncteurs pr\'ec\'edents sont des
isomorphismes de cat\'egories fibr\'ees qui commute avec l'action
de $L_i$.

\medskip

On munit ainsi la classe des tours tors\'ees de structures de
cat\'egories, dont les fl\`eches sont les $n-i-$morphismes qu'on
appelle cat\'egorie des $n-i-$champs et qu'on d\'enote par
$Cham(n-i)$.

\medskip

Soit $C$ une cat\'egorie, consid\'erons une sous-cat\'egorie
$Cham(C)$ de la cat\'egorie des champs appel\'ee encore
cat\'egorie des faisceaux de cat\'egories de base $C$,  telle que

(i) Les objets de $Cham(C)$ sont des gerbes de lien commutatif.

(ii) Munie de la structure de cat\'egorie dont les morphismes sont
 les morphismes de gerbes, la cat\'egorie $Cham(C)$ est ab\'elienne.

\medskip

La classification des extensions $0\rightarrow {F'_1}\rightarrow
F_{n-1}..\rightarrow F_1\rightarrow 0$, est un probl\`eme r\'esolu
dans les cat\'egories ab\'eliennes. Les classes d'isomorphismes
sont d\'etermin\'ees par $Ext^{n-1}({F'}_1,F_1))$.

\bigskip

{\bf Les cocycles classifiant d'une tour de torseurs.}

\bigskip

On va associer \`a une tour de torseurs ab\'elienne
$F_n\rightarrow F_{n-1}...F_2\rightarrow F_1$, une suite de
$l-$cocycles $c_l$ \`a valeurs dans $L_l$, dont la classe de
cohomologie
 appartenant \`a $H^l(F_1,L_l)$ ne d\'epend pas de la
classe d'isomorphisme de la tour.

La famille couvrante de la topologie de $F_1$, $(X_i\rightarrow
X)_{i\in I}$, utilis\'ee pour construire le cocycle sera
suppos\'ee contractible, ceci implique que les groupes de
cohomologie de Cech d'un faisceau relatifs \`a ce recouvrement
sont les groupes de cohomologie du faisceau.

 Supposons que la cat\'egorie $F_1$ est un topos, on
d\'esigne par $U$ son objet final. On notera $X_{i_1..i_n}$ le
produit fibr\'e $X_{i_1}\times_U X_{i_2}..\times_U X_{i_n}$, par
${F_l}_{X_{i_1..i_l}}$ les objets de $F_l$ qui se projettent sur
${X_{i_1..i_l}}$ et par la suite de foncteurs $F_l\rightarrow
F_{l-1}..\rightarrow F_1$. Pour un objet $x_{i_1..i_l}$ de
${F_l}_{X_{i_1..i_l}}$, on notera ${x_{i_1..i_l}}^{i_{l+1}..i_j}$
sa restriction \`a $X_{i_1..i_j}$.

\medskip

 Consid\'erons la fl\`eche

$$
g_{ij}:x_j^i\rightarrow x_i^j
$$
et d\'efinissons  la $2-$chaine

$$
c_{i_1i_2i_3}=g_{i_3i_1}g_{i_1i_2}g_{i_2i_3}
$$

 \`a valeurs dans $L_2$. Giraud [Gi] a d\'emontr\'e que $c_{i_1i_2i_3}$
 est un
$2-$cocycle. La chaine $c_{i_1i_2i_3}$ est un morphisme de la
restriction ${x_{i_3}}^{i_1i_2}$, de $x_{i_3}$ \`a
$X_{i_1i_2i_3}$.

\medskip

Soit $x_{i_1i_2i_3i_4}$ un objet de
${{F_3}_{x_{i_3}}}^{i_1i_2i_4}$.

Les morphismes $c_{i_1i_2i_3}$, $c_{i_1i_2i_4}$, $c_{i_1i_3i_4}$,
...,$c_{i_2i_3i_4}$ de ${x_{i_3}}^{i_1i_2i_4}$ se rel\`event en
des morphismes $(c_{i_1i_2i_3})^*$, $(c_{i_1i_2i_4})^*$,
$(c_{i_1i_3i_4})^*$,...,$(c_{i_2i_3i_4})^*$ de $x_{i_1i_2i_3i_4}$
d'apr\`es l'axiome $(4)$.

On peut d\'efinir

$$d(c_2)^* =(c_{i_2i_3i_4})^*-(c_{i_1i_3i_4})^* +
(c_{i_1i_2i_4})^*- (c_{i_1i_2i_3})^* .
$$

 Puisque $c_2$ est une chaine, on d\'eduit de l'axiome 3 qu'il existe un \'el\'ement
  $c_{i_1..i_4}$ appartenant \`a $L_3(X_{i_1}\times_X..\times X_{i_4})$
  tel que $d(c_2)^*=c_{i_1i_2i_3i_4}=c_3(X_{i_1i_2i_3i_4})$.

\bigskip

{\bf Proposition 2.}

{\it La chaine $c_3$ est un cocycle de Cech.}

\medskip

{\bf Preuve.}

On a

$$
\sum_{l=1}^{l=5}(-1)^lc_{i_1..\hat i_l..i_5}
$$

$$
=\sum_{l=1}^{l=5}(-1)^l\sum_{j=1}^{j=5}(-1)^j(c_{i_1..\hat
i_j..\hat i_l..i_5}^*)=0
$$

La derni\`ere \'egalite r\'esulte de l'axiome $(6)$.

\medskip

Supposons d\'efini le cocycle $c_{l}$, repr\'esent\'e par la
chaine $c_{i_1...i_{l+1}}$, automorphisme de $x_{i_1...i_{l+1}}$.
Consid\'erons un objet $x_{i_1...i_{l+2}}$ de la fibre
${F_{l+1}}_{{x_{i_1...i_{l+1}}}}^{i_{l+2}}$. Les morphismes
$c_{i_2..i_{l+2}}$, $c_{i_1i_3..i_{l+2}}$,..,$c_{i_1..i_{l+2}}$ se
rel\`event en des morphismes $(c_{i_2..i_{l+2}})^*$,
$(c_{i_1i_3..i_{l+2}})^*$ ,..,$(c_{i_1..i_{l+1}})^*$ de
$x_{i_1..i_{l+2}}$, d'apr\`es l'axiome $(3).$

$$
\sum_{j=1}^{j=l+2}(-1)^j(c_{i_1..\hat i_j..i_{l+2}})^*=
c_{i_1..i_{l+2}}
$$

o\`u $c_{i_1..i_{l+2}}$  est un \'el\'ement de
$L_{l+1}(X_1\times_X..\times_XX_{l+2})$.  Cette derni\`ere
affirmation r\'esulte du fait que la chaine $c_{i_1..i_{l+1}}$ est
un cocycle, et de l'axiome $3$.

\bigskip

{\bf Proposition 3.}

{\it La chaine $c_{i_1...i_{l+2}}$ est un cocycle.}

\medskip

{\bf Preuve.}

On a:

$$
\sum_{j=1}^{j=l+3}(-1)^jc_{i_1..\hat i_j..i_{l+3}}
$$

$$
=\sum_{j=1}^{j=l+3}(-1)^l\sum_{d=1}^{d=l+3}(-1)^jc_{i_1..\hat
i_j..\hat i_l..i_{l+3}}^*=0.
$$

\medskip

La derni\`ere \'egalit\'e r\'esulte de l'axiome $(6)$.

\medskip

On va maintenant montrer que les classes de cohomologie des
cocycles $c_2,..,c_{n}$ ne d\'ependent pas des diff\'erents choix
effectu\'es pour les construire.

\medskip

{\bf Proposition 4.}

{\it La classe de cohomologie  du cocycle $c_n$ est ind\'ependant
des diff\'erents choix effectu\'es   pour la d\'efinir.}

\medskip

{\bf Preuve.}

Les r\'esultats de Giraud [Gi] prouvent que le cocycle $c_2$ ne
d\'epend pas des diff\'erents choix effectu\'es pour le
construire. Supposons qu'il en ait de m\^eme pour le cocycle
$c_l$.

Montrons que $c_{l+1}$ ne d\'epend pas des rel\`evements des
$(c_{i_1..i_{l+1}})^*$ effectu\'es, les $x_{i_1..i_{l+2}}$ \'etant
choisis. Consid\'erons d'autres rel\`evement
$(c'_{i_1..i_{l+1}})^*$ de $c_{i_1..i_{l+1}}$. Le morphisme
$(c'_{i_1..i_{l+1}})^* - (c_{i_1..i_{l+1}})^*$ est un rel\`evement
de l'identit\'e de $x_{i_1..i_{l+1}}$. On en d\'eduit de l'axiome
$3$, l'existence d'un \'el\'ement $h_{i_1..i_{l+1}}$ de
$L_{l+1}(X_{i_1..i_{l+1}})$ tel que
$(c'_{i_1..i_{l+1}})^*=(c_{i_1..i_{l+1}})^* + h_{i_1..i_{l+1}}$.
Il en r\'esulte que la diff\'erence entre les cocycles d\'efinis
par les rel\`evements $(c_{i_1..i_{l+1}})^*$ et
$(c'_{i_1..i_{l+1}})^*$ est un bord.

Montrons maintenant que le cocycle ne d\'epend pas du choix des
$x_{i_1..i_{l+2}}$. Consid\'erons un autre choix d'\'el\'ements
$x'_{i_1..i_{l+2}}$ dans ${F_{l+1}}_{x"_{i_1..i_{l+1}}}$. On note
respectivement $(c_{i_1..i_{l+1}})^*$ et $(c'_{i_1..i_{l+1}})^*$,
les rel\`evements de $c_{i_1..i_{l+1}}$ \`a ${x_{i_1..i_{l+2}}}$
et $x'_{i_1..i_{l+2}}$.

Consid\'erons un morphisme  $u:x_{i_1..i_{l+2}}\rightarrow
x'_{i_1..i_{l+2}}$.

$$
{u}(\sum_{j=1}^{j=l+2}(-1)^j{(c_{i_1..\hat
i_j..i_{l+2}})^*})u^{-1}
$$

$$
=\sum_{j=1}^{j=l+2}(-1)^j{u}{(c_{i_1..\hat i_j..i_{l+2}}^*})u^{-1}
$$

$$
(\sum_{j=1}^{j=l+2}(-1)^j(c'_{i_1..\hat i_j..i_{l+2}})^*+ d(c)
$$
 o\`u $d(c)$ est le bord d'une chaine.

$$
u(\sum_{j=1}^{j=l+2}(-1)^j{(c_{i_1..\hat i_j..i_{l+2}})^*})u^{-1}
=c_{i_1..i_{l+2}}
$$

d'apr\`es l'axiome $3$. On en d\'eduit  le r\'esultat.

\bigskip

 {\bf Proposition 5.}

{\it  Les classes de cohomologie des cocycles associ\'es \`a deux
tours de torseurs ab\'eliennes isomorphes et de m\^eme liens
$L_2,..,L_n$  sont identiques.}

\medskip

{\bf Preuve.}

Soient $ u={F_n}\rightarrow {F_{n-1}}...\rightarrow {F_1}$, et
$u'={F'_n}\rightarrow {F'_{n-1}}...\rightarrow {F'_1}$, deux tours
de torseurs isomorphes, par le biais de la famille d'isomorphismes
fonctoriels $f_i:{F_i}\rightarrow {F'_i}$. On d\'enote
respectivement par $(c_2,...,c_{n})$ et $(c'_2,...,c'_{n})$ les
familles de cocycles associ\'es aux tours $u$ et $u'$. D'apr\`es
la th\'eorie de Giraud, les classes de cohomologie $[c_2]$ et
$[c'_2]$ sont identiques. Supposons que les classes de cohomologie
$[c_i]$ et $[c'_i]$ sont identiques et montrons que les classes
$[c_{i+1}]$ et $[c'_{i+1}]$ coincident.

La famille $(u_1(X_i)\rightarrow u_1(X))_{i\in I}$ est une famille
g\'en\'eratrice contractible de la topologie de $F'_1$.
L'isomorphisme $f_{i+1}:F_{i+1}\rightarrow F'_{i+1}$ induit un
isomorphisme
${f_i}_{x_{i_1..i_l}}:{F_{i+1}}_{x_{i_1..i_{l+1}}}\rightarrow
{F'_{i+1}}_{x_{i_1..i_{l+1}}}$. La famille
${f_i}_{x_{i_1..i_l}}(x_{i_1..i_{l+2}})=x'_{i_1..i_{l+2}}$ peut
\^etre choisie pour d\'efinir le cocycle $c'_{i_1..i_{l+2}}$
d'apr\`es la proposition $3$, cette m\^eme proposition entraine
que
${f_i}_{x_{i_1..i_l}}(c_{i_1..i_{l+2}})^*{f_i}_{x_{i_1..i_l}}^{-1}$
est un rel\`evement de $c'_{i_1..i_{l+1}}$. On en d\'eduit le
r\'esultat.

\bigskip

Nous allons maintenant \'enoncer l'axiome de recollement des
objets d'une tour tors\'ee, nous avons besoin dans cette optique
de la d\'efinition suivante:

\medskip

{\bf D\'efinition 8.}

Soit $ u={F_n}\rightarrow {F_{n-1}}...\rightarrow {F_1}$, une tour
de torseurs. On dira que la tour $u$ est triviale s'il existe une
tour de torseurs $u'=F'_{n-1}\rightarrow {F'_{n-2}}...\rightarrow
{F'_1}$, une application croissante

$$
j:\{1,...,n-1\}\longrightarrow \{1,...,n\}
$$

telle que $j(n-1)=n$, des morphismes

$$
f_i:F'_i\longrightarrow F_{j(i)}
$$

tels que le diagramme suivant soit commutatif:

$$
\matrix{ F'_i{f_i\over\longrightarrow} F_{j(i)}\cr \downarrow &&
\downarrow \cr F'_{i-1} {f_{i-1}\over\longrightarrow} F_{j(i-1)}}
$$

\medskip

{\bf Axiome de recollement des objets.}

\medskip

La tour de torseurs est triviale si et seulement si la classe de
cohomologie du cocycle $[c_{n}]$ est nulle.

\medskip

Consid\'erons la cat\'egorie fibr\'ee $F_{i+1}\rightarrow F_i$,
d'apr\`es l'axiome $(6)$ il existe un faisceau ab\'elien
$L'_{i+1}$ sur $F_1$ tel que $Aut({F_{i+1}}_{x_i})$ l'ensemble des
morphismes de ${F_{i+1}}_{x_i}$ ne se projettant pas  forc\'ement
sur l'identit\'e est l'ensemble des sections
$L'_{i+1}(p_2..p_i(x_i))$ de $L'_{i+1}$. On a alors la suite
exacte

$$
0\rightarrow L_{i+1}\rightarrow L'_{i+1}\rightarrow
L_{i}\rightarrow 0
$$
On en d\'eduit la suite exacte en cohomologie:

$$
...\rightarrow H^l(F_1,L_{i+1})\rightarrow
H^l(F_1,L'_{i+1})\rightarrow H^{l}(F_1,L_{i})\rightarrow
H^{l+1}(F_1,L_{i+1})...
$$

\medskip

{\bf Proposition 9.}

{\it Le cocycle $c_{l+1}$ est l'image de $c_{l}$ par l'op\'erateur
de Dolbeault.}

\medskip

{\bf Preuve.}

La construction de la classe de cohomologie $c_{l+1}$  est celle
de l'image de $c_{l}$ par l'op\'erateur de Dolbeault.

\bigskip

Une question importante est de d\'eterminer les classes de
cohomologie de $H^ n(F_1,L_n)$ qui peuvent \^etre r\'ealis\'ees
comme classes classifiantes d'une tour tors\'ee. La proposition
pr\'ec\'edente  donne une interpr\'etation g\'eom\'etrique de
l'op\'erateur de Dolbeault et une condition n\'ecessaire \`a cette
r\'ealisation. Cette condition est en fait aussi suffisante.

\medskip

{\bf Th\'eor\`eme  2.}

{\it Soient $F_n\rightarrow F_{n-1}...\rightarrow F_1$ une tour de
torseurs associ\'ee aux faisceaux $L_2,...,L_n$, et  $L_{n+1}$ un
faisceau d\'efini sur $F_1$ et $c_{n+1}$ un \'el\'ement de
$H^{n+1}(F_1,L_{n+1})$.  On suppose qu'il existe faisceau
$L'_{n+1}$ sur $F_1$ et une suite exacte
$$
0\rightarrow L_{n+1}\rightarrow L'_{n+1}\rightarrow L_n\rightarrow
0
$$
telle que $c_{n+1}$ soit l'image de $c_{n}$ par l'op\'erateur de
Dolbeault. Alors on peut prolonger la tour de torseurs en une tour
de torseurs $F_{n+1}\rightarrow F_n...\rightarrow F_1.$}

\medskip

{\bf Preuve.}

On va construire la cat\'egorie $F_{n+1}$ comme suit:
consid\'erons l'objet $x_{i_1..i_{n+1}}$ de $F_n$, au-dessus de
l'objet $x_{i_1..i_{n}}$ utilis\'e pour construire le cocycle
$c_{n}$. Les objects de $F_{n+1}$ sont les objets des gerbes
triviales ${L'_{n+1}}_{x_{i_1..i_{n+1}}}$ dont chaque section est
isomorphe \`a $L'_{n+1}(X_{i_1..i_{n+1}})$, le foncteur
$F_{n+1}\rightarrow F_n$ associe \`a un \'el\'ement de
${L'_{n+1}}_{x_{i_1..i_{n+1}}}$ l'objet $x_{i_1..i_{n+1}}$. Soient
$x'_{i_1..i_{n+1}}$ un objet de $F_n$, $x_{i_1..i_{n+2}}$ et
$x'_{i_1..i_{n+2}}$ des objets de $F_{n+1}$ respectivement
au-dessus de $x_{i_1..i_{n+1}}$, et $x'_{i_1..i_{n+1}}$.
L'ensemble $Hom_{F_{i+1}}(x_{i_1..i_{n+2}},x'_{i_1..i_{n+2}})$ se
d\'efinit comme suit:
$Hom_{F_{n+1}}(x_{i_1..i_{n+2}},x'_{i_1..i_{n+2}})$, est
l'ensemble $L'_{n+1}(X_{i_1..i_{n+1}})_u$ o\`u $u$ est une
fl\`eche entre $x_{i_1..i_{n+1}}$ et $x'_{i_1..i_{n+1}}$. On note
$(l,u)$ un de ses \'el\'ements, $l\in L'_{n+1}(X_{i_1..i_{n+1}})$
$(l,u)\circ (l',u')=(l+l',uu')$.

\bigskip

{\bf Proposition 10.}

{\it Soit $(c_2,...,c_n)$ la famille de cocycles classifiant d'une
tour tors\'ee $F_n\rightarrow F_{n-1}\rightarrow...\rightarrow
F_1$. La nullit\'e de la classe de cohomologie du cocycle $c_l$
entraine celles des cocycles $c_j, j\geq l$.}

\medskip

{\bf Preuve.}

Supposons que le cocycle $c_l$ est d\'efini par les chaines
$c_{i_1..i_{l+1}}$, et que sa classe de cohomologie soit nulle. On
peut alors supposer que la chaine $c_{i_1..i_{l+1}}$ est nulle et
choisir $(c_{i_1..i_{l+1}})^*=0$.

\bigskip

On va maintenant appliquer la notion de tour tors\'ee \`a notre
probl\`eme initial. On consid\`ere donc un fibr\'e principal de
groupe structural $L_0$, de base la vari\'et\'e diff\'erentielle
$N$, et une suite d'extensions \`a noyaux commutatifs.

$$
1\rightarrow H_1\rightarrow L_1\rightarrow L_0\rightarrow 1
$$
...
$$
1\rightarrow H_n\rightarrow L_n\rightarrow L_{n-1}\rightarrow 1.
$$

telles que le quotient $L_{i+1}/L_{i-1}$ soit commutatif. et
l'application $L_{i+1}\rightarrow L_i$ a des sections locales.

Soit $(U_i)_{i\in I}$ un recouvrement de $N$ par des ouverts
contractibles. On d\'efinit une tour de torseurs comme suit:

La cat\'egorie $F_0$ est celle des ouverts de $N$,

Les objets de la cat\'egories $F_1$ sont les fibr\'es principaux
de groupe structural $L_1$ au-dessus des ouverts de $N$. Pour tout
ouvert $U$ de $N$, on d\'enote par ${F_1}_U$ les objets de $F_1$
de base $U$. On suppose que le quotient des objets de ${F_1}_U$
par $H_1$ est la restriction du $L_0-$fibr\'e principal \`a $U$.
Les morphismes entre objets de ${F_1}_U$ sont les isomorphismes de
fibr\'es qui se projettent sur l'identit\'e. Les objets de
${F_1}_{U_i}$ sont isomorphes \`a $U_i\times L_1$. On d\'efinit
ainsi une gerbe de lien le faisceau des fonctions de $N$ \`a
valeurs dans $H_1$.

\medskip

Supposons d\'efinie la cat\'egorie $F_i$. Un  objet $e_{i+1}$ de
$F_{i+1}$ est un fibr\'e principal de groupe structural $L_{i+1}$
de base un objet $e_i$ de $F_i$, tel que le quotient de $e_{i+1}$
par $H_{i+1}$ est $e_i$. Soient $e_{i+1}$ et $e'_{i+1}$ deux
objets de $F_{i+1}$ qui se projettent sur $e_i$. L'ensemble
$Hom(e_{i+1},e'_{i+1})$ est constitu\'e des morphismes de fibr\'es
qui se projettent sur l'identit\'e de $e_i$. Les axiomes qui
d\'efinissent la notion de tour de torseurs sont satisfaits par la
suite de foncteurs que nous venons de d\'efinir:

L'axiome $(1)$ est satisfait car la cat\'egorie $F_1$ est une
gerbe de base la cat\'egorie des ouverts de $N$.

L'axiome $2$ est satisfait car les fl\`eches des cat\'egories
$F_i$ sont des morphismes de fibr\'es inversibles.

L'axiome $6$ est satisfait car on a suppos\'e que
$L_{i+1}/L_{i-1}$ est commutatif.

Il est \'evident de v\'erifier le reste des axiomes.

\medskip

Le cocycle $c_{i+1}$ est \`a valeurs dans le faisceau des
fonctions d\'efinies sur $N$ et \`a valeurs dans $L_i$.

\medskip

Plus g\'en\'eralement on peut consid\'erer les tour de gerbes
$F_n\rightarrow F_{n-1}\rightarrow...F_1\rightarrow F_0$, o\`u les
fibres de $F_{1}\rightarrow F_0$ sont des gerbes triviales dont
les objets sont isomorphes \`a des fibr\'es principaux au-dessus
des ouverts de $N$, supposons que la cat\'egorie $F_i$ est
d\'efinie. La fibre de $e_i$ de $F_i$ pour le foncteur
$F_{i+1}\rightarrow F_i$ est une gerbe triviale dont les objets
sont des fibr\'es principaux de groupes structural $L_{i+1}$
au-dessus de $e_i$. Les morphismes \'etant les applications de
fibr\'es qui se projettent sur l'identit\'e. Ces tours ne sont pas
n\'ecessairement d\'efinies par un probl\`eme de rel\`evement.

\bigskip

{\bf II Geom\'etrie diff\'erentielle des tours de torseurs
diff\'erentielles.}

\bigskip

Le but de cette partie est d'\'etudier la g\'eom\'etrie
diff\'erentielle des tours de torseurs, en commencant par
rappeller celle des gerbes d\'efinie par Brylinski. On notera
l'alg\`ebre de Lie d'un groupe $L$ par ${\cal L}$.

\medskip

{\bf G\'eom\'etrie diff\'erentielle des gerbes.}

\medskip

Consid\'erons une gerbe principale $C$ d\'efinie par un fibr\'e
principal $h$ de groupe structural $L$, de base la vari\'et\'e
$N$, et une extension $1\rightarrow H_1\rightarrow L_1\rightarrow
L\rightarrow 1$.

\medskip

{\bf D\'efinition 1.}

 Soit $\omega$ une connexion d\'efinie sur $h$, une structure
 connective sur $C$ est d\'efinie par:

 Pour tout ouvert $U$ de $N$ et tout objet $e_U$ de $C(U)$,
  l'ensemble des connections $Co(e_U)$ image inverse de
 la restriction de $\omega$ \`a $U$. Soit $(V_i)_{i\in I}$ un recouvrement de
 $U$ par des ouverts contractibles. La restriction de la connection $\omega$ est d\'efinies
 par des $1-$formes locales $\omega_i:V_i\rightarrow {\cal L}$ v\'erifiant
 $\omega_j-\omega_i={g_{ij}}^{-1}d(g_{ij})$. Un \'el\'ement $w$ de $Co(e_U)$
est d\'efini par des formes $w_i: V_i\rightarrow {\cal L}_1$
d\'efinissant une connection sur $e_U$ telles que
$\omega_i=h_1\circ w_i$, o\`u $h_1:{\cal L}_1\rightarrow {\cal L}$
est la projection canonique.

 \bigskip

 {\bf Le cocycle caract\'eristique d'une structure connective d'une gerbe.}

 \bigskip

 Consid\'erons un recouvrement de $N$, $(U_i)_{i\in I}$ par des ouverts
 contractibles. Pour tout objet $e_i$ de $C(U_i)$, on choisit un
 \'el\'ement $w_i$ de $Co(e_i)$. Soit
 $g_{ij}:e_j^i\rightarrow e_i^j$ une fl\`eche. La forme
 $\alpha_{ij}= w_j-{g_{ij}}^*w_i$ d\'efinie sur $U_i\cap
 U_j$, est \`a valeurs dans ${\cal H}_1$.

 \medskip

 {\bf Proposition 1.}

 {\it Les formes $\alpha_{ij}$ v\'erifient la relation:
 $$
 \alpha_{jk}-\alpha_{ik}+\alpha_{ij}=c_{ijk}^{-1}dc_{ijk}.
 $$
o\`u $c_{ijk}$ est le cocycle classifiant de la gerbe.}

\medskip

{\bf Preuve.}

On a:

$$
\alpha_{jk}-\alpha_{ik}+\alpha_{ij}=w_k-{g_{jk}}^*w_j -
(w_k-{g_{ik}}^*w_i) + {g_{jk}}^*(w_j-{g_{ij}}^*w_i)
$$

$$
= {g_{ik}}^*(w_i-{{g_{ik}}^*}^{-1}{g_{jk}}^*{g_{ij}}^*w_i)
$$

$$
=w_i - {c_{ijk}}^*w_i={c_{ijk}}^{-1}dc_{ijk}
$$

\bigskip

{\bf Courbure d'une structure connective.}

\medskip

 Pour tout ouvert $U_i$ d'un recouvrement contractible de $N$, et $e_i$
 un objet de $C(U_i)$, consid\'erons un \'el\'ement $w_i$ de
 $Co(e_i)$, la courbure de $w_i$ est la $2-$forme:

 $$
 L(w_i)=d(w_i) + {1\over 2}[w_i,w_i]
 $$

 Tout \'el\'ement de $Co(e_i)$ $w'_i$ s'\'ecrit
 $w'_i=w_i +u_i$, o\`u $u_i$ est une $1-$forme d\'efinie sur $U_i$ \`a
 valeurs dans ${\cal H}_1$, on a

 $$
 L(w'_i)= L(w_i) + d(u_i)
 $$
car $H_1$ est commutatif. Ceci implique que $dL$ est une $3-$forme
d\'efinie sur $N$ \`a valeurs dans ${\cal H}_1$ appel\'ee la
courbure.

\bigskip

{\bf Le cas g\'en\'eral.}

\bigskip

On \'etudie  la g\'eom\'etrie diff\'erentielle d'une tour de
torseurs associ\'ee au probl\`eme d'extension du groupe structural
des fibr\'es principaux.

On consid\`ere donc un fibr\'e principal de groupe structural
$L_0$ de base la vari\'et\'e diff\'erentielle $N$, et une suite
d'extensions centrales:

$$
1\rightarrow H_1\rightarrow L_1\rightarrow L_0\rightarrow 1
$$
...
$$
1\rightarrow H_n\rightarrow L_n\rightarrow L_{n-1}\rightarrow 1.
$$
On associ\'e \`a cette suite une tour de torseurs $F_n\rightarrow
F_{n-1}\rightarrow...\rightarrow F_1\rightarrow F$.

\medskip

{\bf D\'efinition 2.}

Une structure connective sur la tour de torseurs pr\'ec\'edente
est d\'efinie comme suit:

On consid\`ere une connexion $w^0$ sur le $L_0-$fibr\'e principal.

Soit $U$ un ouvert de $N$, pour tout objet ${e^1}_U$ de ${F_1}_U$,
on consid\`ere l'ensemble $Co^1({e^1}_U)$ des connexions
d\'efinies sur ${e^1}_U$ qui se projettent sur la restriction de
$w_0$.

Supposons d\'efini l'ensemble de connexions $Co^i(e_i)$ pour tout
objet $e_i$ de $F_i$. Soit $e_{i+1}$ un objet de $F_{i+1}$
au-dessus de $e_i$, les \'el\'ements de $Co^{i+1}(e_{i+1})$ sont
les connections d\'efinient sur $e_{i+1}$ qui se projettent sur un
\'el\'ement de $Co^i(e_i)$.

\bigskip

{\bf Courbure d'une tour de torseurs.}

\bigskip

Soit $F_n\rightarrow F_{n-1}\rightarrow..F_1\rightarrow F$ une
tour de torseurs d\'efinie par un probl\`eme d'extension.

Pour tout objet $e_i$ de ${F_1}_{U_i}$, on consid\`ere une
connection $w_i$, et $L(w_i)$ sa courbure. On a
$$
d(L(w_i))=d(L(w_j)),
$$

ceci permet de d\'efinir la trois forme $L^3$ dont la restriction
\`a $U_i$ est $d(L(w_i))$. La forme $L^3$ est la courbure de la
gerbe $F_1\rightarrow F$.

Supposons d\'efinie la courbure de la tour
$F_j\rightarrow..\rightarrow F_1\rightarrow F$, c'est une
$j+2-$forme $\Omega^{j+2}$ \`a valeurs dans ${\cal L}_j$.

On consid\`ere un recouvrement $(U_i)_{i\in I}$ de $N$ par des
ouverts contractibles. Sur chaque $U_i$, on consid\`ere une
$j+2-$forme $c_i$ \`a valeurs dans ${\cal L}_{j+1}$ au-dessus de
la restriction de $\Omega^{j+2}$ sur $U_i$.

La chaine $c_{i_2}-c_{i_1}$ est un $1-$cocycle de Cech \`a valeurs
dans l'ensemble des $j+2-$formes d\'efinies sur $N$ et \`a valeurs
dans ${\cal H}_{j+1}$.

Le morphisme de De Rham-Cech l'identifie \`a une $j+3-$forme qui
est la courbure de la tour de torseurs
$F_{j+1}\rightarrow...\rightarrow F$.

On d\'efinit ainsi la courbure $L$ de la tour de torseurs qui est
une $n+2-$forme \`a valeurs dans ${\cal H}_{n}$.

\bigskip

{\bf Une chaine diff\'erentielle fondamentale de la tour de
torseurs.}

\bigskip

 Soient  $c_{i_1..i_{n+2}}$ le $n+1-$cocycle classifiant de la tour de torseurs
 $F_n\rightarrow F_{n-1}...\rightarrow F_1\rightarrow F$. On pose
 $u_{i_1..i_{n+2}}=log(c_{i_1..i_{n+2}})$. La chaine $u_{i_1..i_{n+2}}$ est un
 $n+1-$cocycle \`a valeurs dans ${\cal H}_n$. On d\'eduit du th\'eor\`eme
 de De Rham-Cech l'existence d'un $n-$cocycle de $1-$formes $u^1_{i_1..i_{n+1}}$
 tel que

 $$
 d(u_{i_1..i_{n+2}})=\sum_{j=1}^{j=n+2}(-1)^j{u^1}_{i_1..\hat i_j..i_{n+2}}
 $$

Supposons d\'efinis les $l-$formes $u^l_{i_1..i_{n-l+2}}$ tels que

$$
d({u^{l-1}}_{i_1..i_{n-l+3}})=\sum_{j=1}^{j=n-l+3}(-1)^j{u^l}_{i_1..\hat
i_j..i_{n-l+3}}
 $$

On en d\'eduit que

$$
\sum_{j=1}^{j=n-l+3}(-1)^jd({u^l}_{i_1..\hat i_j..i_{n-l+3}})=0
 $$

On en d\'eduit du th\'eor\`eme de Poincare-Cech l'existence de
$l+1-$formes ${u^{l+1}}_{i_1..i_{n-l+1}}$ telles que

$$
d({u^{l}}_{i_1..i_{n-l+2}})=\sum_{j=1}^{j=n-l+2}(-1)^j{u^{l+1}}_{i_1..\hat
i_j..i_{n+2}}
 $$

On en d\'eduit par r\'ecurence l'existence d'une $n+2-$forme
$u^{n+2}$ forme sur $N$.

\medskip

{\bf Proposition 2.}

{\it La classe de cohomologie de la courbure de la tour de
torseurs coincide avec celle de la forme $u^{n+2}$  d\'efinie au
paragraphe pr\'ec\'edent.}

\medskip

{\bf Preuve.}

La preuve se fait par r\'ecurence. Consid\'erons la tour de
torseurs $F_n\rightarrow F_{n-1}..\rightarrow F_1\rightarrow F$,
associ\'ee \`a un probl\`eme d'extension. Rapellons que la
courbure de la gerbe $F_1\rightarrow F$ est d\'efinie en
consid\'erons un rev\^etement contractible $(U_i)_{i\in I}$, un
objet $e_i$ de ${F_1}_{U_i}$, et une connection $w_i$  de
$Co(e_i)$ La courbure $L(w_i)$ de $w_i$ v\'erifie $L(w_j)=L(w_i) +
d(u_{ij})$, o\`u $u_{ij}=w_j-{g_{ij}}^*w_i$, de plus on a la
relation $u_{i_2i_3}-u_{i_1i_3}+
u_{i_1i_2}=(c_{i_1i_2i_3})^{-1}d(c_{i_1i_2i_3})$. La courbure est
la forme donc la restriction \`a $U_i$ est $dL(w_i)$. Ceci
signifie que la courbure de la tour est la forme du paragraphe
pr\'ec\'edent associ\'ee \`a la gerbe $F_1\rightarrow F$.

Supposons que l'\'enonc\'e v\'erifi\'e pour la tour de torseurs
$T_{n-1}=F_{n-1}\rightarrow F_{n-2}...\rightarrow F_1\rightarrow
F$. La construction de la courbure $L$ initiale prouve que sa
classe de cohomologie est l'image de celle de la courbure $L'$ de
$T_{n-1}$ par l'homomorphisme $H^{n}(N,L_{n-1})\rightarrow
H^{n+1}(N,L_n)$.

Les constructions des formes $u^{n+1}$ et $u^{n+2}$ associ\'ees
aux tours $T_{n-1}$ et $T_n$ montrent que $[u^{n+1}]$ est l'image
de $[u^{n}]$ par le m\^eme homomorphisme pr\'ec\'edent. On en
d\'eduit le r\'esultat.

\bigskip

{\bf Classes charact\'eristiques des tours de torseurs.}

\bigskip

Soit $T_n=F_n\rightarrow F_{n-1}..\rightarrow F_1\rightarrow F$
une tour tors\'ee d\'efinie par un probl\`eme d'extension,  munie
d'une structure connective.

La courbure de $T_n$ est une $n+2-$forme \`a valeurs dans ${\cal
H}_n$. Pour tout polynome de ${\cal H}_n$, $P$, on d\'efinit la
forme $P(L)$ qui est une forme caract\'eristique de la tour de
torseur.

\bigskip

{\bf Holonomie d'une tour de torseurs.}

\bigskip

Soit $T_n=F_n\rightarrow F_{n-1}..\rightarrow F_1\rightarrow F$
une tour de torseurs associ\'ee \`a un probl\`eme d'extension,
munie d'une structure connective. Rappellons que la courbure est
d\'efinie comme suit:

On consid\`ere  $c_{i_1..i_{n+2}}$ le $n+1-$cocycle classifiant de
la tour de torseurs
 $F_n\rightarrow F_{n-1}...\rightarrow F_1\rightarrow F$. On pose
 $u_{i_1..i_{n+2}}=log(c_{i_1..i_{n+2}})$. Il existe des $l-$formes $u^l_{i_1..i_{n-l+2}}$
 telles que

$$
d({u^{l-1}}_{i_1..i_{n-l+3}})=\sum_{j=1}^{j=n-l+3}(-1)^j{u^l}_{i_1..\hat
i_j..i_{n-l+3}}
 $$

Supposons que la courbure $u^{n+2}$ soit nulle, ceci implique que

$$
{u^{n+2}}_i=d({v^{n+1}}_i)
$$

$$
{u^{n+1}}_{i_1i_2}=
{v^{n+1}}_{i_2}-{v^{n+1}}_{i_1}+d({v^{n}}_{i_1i_2})
$$
...
$$
u^l_{i_1...i_{n-l+2}}=\delta (v^l_{i_1..i_{n-l+1}})+
d({v^{l-1}}_{i_1..i_{n-l+2}})
$$
...
$$
u^1_{i_1...i_{n+1}}=\delta({v^1}_{i_1..i_{n}}) +
d({v}_{i_1...i_{n+1}})
$$

\medskip

On d\'enote par $w_{i_1..i_{n+2}}$ la classe
$(c_{i_1..i_{n+2}})^{-1}\delta({v}_{i_1...i_{n+1}})$. C'est la
classe d'holonomie de la tour.

\medskip

Soit $V$ une vari\'et\'e de dimension $n$ sans bord et
$h:V\rightarrow N$ une application diff\'erentiable. La courbure
de la tour de torseur $h^*(T_n)$ est nulle.  On peut d\'efinir le
cocycle d'holonomie $h_{i_1..i_{n+2}}$. La dimension de $V$
entraine que le cocycle $log(h_{i_1..i_{n+2}})$ est trivial, il
est donc le bord d'un cocycle $h'_{i_1..i_{n+1}}$, tel que
$d(h'_{i_1..i_{n+1}})=0$. On peut alors utiliser l'isomorphisme de
De Rham Cech et identifier $h'_{i_1..i_{n+1}}$ \`a une $n-$forme
$Hol(V)$ de $V$ \`a valeurs dans ${\cal H}_n$.
 L'holonomie autour de $V$ est:

$$
\int_V Hol(V)
$$

\bigskip

{\bf Applications \`a la physique th\'eorique.}

\bigskip

Soit $x$ une particule  se mouvant sur une vari\'et\'e
riemmanienne $(N,<,>)$, lorsque la particule est libre, sa
trajectoire est celle des g\'eod\'esiques de $<,>$. Supposons que
$x$ soit soumis \`a un champ de forces $L$, repr\'esent\'e par une
$2-$forme ferm\'ee de $N$ qu'on note encore $L$. Si la deux forme
$L=d(V)$ l'action de la formulation variationnelle du mouvement de
$x$ est de la forme:

$$
\int_NV=0
$$
La forme $L$ n'est pas toujours exacte, le mouvement de $x$ \`a
tout de m\^eme une solution variationnelle si la classe deux
cohomologie de $L$ est enti\`ere. Dans ce cas $L$ est la
premi\`ere classe de Chern d'un fibr\'e en cercle au-dessus de
$N$. Il existe une connection $\Omega$ sur ce fibr\'e d\'efini sur
un recouvrement contractible $(U_i)_{i\in I}$ de $N$ par des
$1-$formes $w_i$, l'action de $x$  sur un chemin $u$ est de la
forme

$$
\int_u u^*Hol(w_i,L)
$$
o\`u $Hol(w_i,L)$ est l'holonomie de la connection.

\medskip

En physique th\'eorique les cordes ont \'et\'e remplac\'ees les
particules, une corde est une application de l'intervalle
$c:I\rightarrow N$, les trajectoires des cordes sont des surfaces
immers\'ees dans $N$. L'action d'une corde  dans le mod\`ele $WZW$
est de la forme

$$
\int u^*Hol(w_i,u_{ij},L)
$$

o\`u $Hol(w_i,u_{ij},L)$ est l'holonomie d'une gerbe le long d'une
surface plong\'ee dans $N$ par l'homomorphisme $u$. Cette action
est celle de $WZW$.

 \medskip

 Plus g\'en\'eralement en th\'eorie des  branes, on \'etudie les mouvements des
 immersions de $V\rightarrow N$, qui repr\'esentent les conditions limites de
 l'\'evolution des cordes. La formulation variationnelle de ce probl\`eme
 est donn\'ee par une action g\'en\'eralisant celle de $WZW$ par

 $$
 \int_V Hol
 $$

o\`u est l'holonomie d'une tour de torseus associ\'ee \`a un probl\`eme de rel\`evement.

\bigskip

\bigskip

\centerline{\bf Bibliographie.}

\bigskip

[Bre] Breen, L. On the classification of $2-$gerbes and
$2-$stacks. Asterisque, 225 1994.

[Br] Bredon, G. E. Sheaf theory. McGraw-HillBook Co., 1967.

[Bry] Brylinski, J.L Loops spaces, Characteristic Classes and
Geometric Quantization, Progr. Math. 107, Birkhauser, 1993.

[Br-Mc] Brylinski, J.L, Mc Laughlin D.A, The geometry of degree
four characteristic classes and of line bundles on loop spaces I.
Duke Math. Journal. 75 (1994) 603-637.

[Ca-Mi] Alan L. Carey, Jouko Mickelsson The universal gerbe,
Dixmier-Douady class, and gauge theory
     Lett.Math.Phys. 59 (2002) 47-60

[De] Deligne, P. Theorie de Hodge III, Inst. Hautes Etudes Sci.
 Math. 44 (1974), $5-77$.

[Du] Duskin, J. An outline of a theory of higher dimensional
descent, Bull. Soc. Math. Bel. S\'erie A 41 (1989) 249-277.

[Fri] Fried, D. Closed similarity affine manifolds. Comment. Math.
Helv. 55 (1980) 576-582.

[Ga-Reis] Gawedzki, Reis N. WZW branes and gerbes Rev.Math.Phys.
14 (2002) 1281-1334

[Gi] Giraud, J. Cohomologie non ab\'elienne.

[God] Godement R. Topologie alg\'ebrique et th\'eorie des
faisceaux. (1958) Hermann.

Grothendieck, Seminaire de Geometrie Algebrique.

[J] Johnson. D-branes. Cambridge University Press

[Mc] Maclane, S. Homology. Springer-Verlag, 1963.

 [Mu-S] M. K. Murray, D. Stevenson (University of Adelaide)
     Commun.Math.Phys. 243 (2003) 541-555

[P] Polchinski, J. String theory. Cambridge University Press

[St] Stuart, J. Constructions with bundle gerbes. Ph D Thesis.
University of Adelaide

[T1] Tsemo, A. Non abelian cohomology: the point of view of gerbed tower

[W] Witten, E. Quantum field theory and the Jones Polynomial,
Comment. Math. Phys. 121 (1989) 351-399.

[Z] Zinn, J. Quantum field theory and critical phenomea. Oxford
Sciences Publications

\end{document}